\documentclass[10pt]{article}
\usepackage{mathrsfs}
\usepackage{epic,eepic,epsf,epsfig}
\usepackage{amsfonts,srcltx,mathrsfs}
\usepackage{amsmath}
\usepackage{amssymb}
\usepackage{amsbsy}
\usepackage{graphicx}
\usepackage{setspace}
\usepackage{amsfonts}
\usepackage{color}
\usepackage{array}

\newtheorem{lem}{Lemma}[section]%
\newtheorem{theorem}[lem]{Theorem}%
\newtheorem{defi}[lem]{Definition}%

 \def\b{\beta}

\def\nd{\mathrel{\bigm|\kern-.7em/}}

\def\f{\noindent}

\newcommand{\qed}{\mbox{\raisebox{0.7ex}{\fbox{}}} \vspace{4truemm}}

\begin{document}
\begin{spacing}{1.0}

\markboth{M-M Gu et. al}{Fault-tolerance of balanced hypercubes with faulty vertices and faulty edges}

\title{Fault-tolerance of balanced hypercubes with faulty vertices and faulty edges}

\author{ \\ Mei-Mei Gu\\
{\small\em Department of Mathematics, Beijing Jiaotong University}, \\{\small\em Beijing 100044, China}\\ {\small\em gum2012@bjtu.edu.cn}\\[0.5cm]
Rong-Xia Hao\footnote{Corresponding author}\\
{\small\em Department of Mathematics, Beijing Jiaotong University},\\ {\small\em Beijing 100044, China}\\ {\small\em rxhao@bjtu.edu.cn}\\}

\footnotetext[1]{This paper was accepted by Ars Combinatoria, March 1, 2016}

\date{}
\maketitle

\f {\bf Abstract}\quad Let $F_{v}$ (resp. $F_{e}$) be the set of faulty vertices (resp. faulty edges) in the $n$-dimensional balanced hypercube $BH_{n}$.
Fault-tolerant Hamiltonian laceability in $BH_{n}$ with at most $2n-2$ faulty edges is obtained in [Inform. Sci. 300 (2015) 20--27].
The existence of edge-Hamiltonian cycles in $BH_{n}-F_{e}$ for $|F_{e}|\leq 2n-2$ are gotten in [Appl. Math. Comput. 244 (2014) 447--456].
Up to now, almost all results about fault-tolerance in $BH_{n}$ with only faulty vertices or only faulty edges.
In this paper, we consider fault-tolerant cycle embedding of $BH_{n}$ with both faulty vertices and faulty edges,
and prove that there exists a fault-free cycle of length $2^{2n}-2|F_{v}|$ in $BH_{n}$ with $|F_{v}|+|F_{e}|\leq 2n-2$ and $|F_{v}|\leq n-1$ for $n\geq 2$.
Since $BH_{n}$ is a bipartite graph with two partite sets of equal size, the cycle of a length $2^{2n}-2|F_{v}|$ is the longest in the worst-case.
\bigskip

\f {\bf Keywords}\quad  Balanced hypercube; Cycle embedding; Fault tolerance; Interconnection network.\\

\section{Introduction}

In a multiprocessor system, processors communicate by exchanging messages through an interconnection network
whose topology often modeled by an undirected graph $G=(V,E)$, where every vertices in $V$ corresponds to a processor, and
every edge in $E$ corresponds to a communication link.

The $n$-dimensional balanced hypercube $BH_n$, proposed
by Wu and Huang \cite{W}, is an important class of generalizations of the popular hypercube interconnection
network for parallel computing. It has many desirable properties, such as bipartite, regularity, recursive structure, vertex-transitive \cite{H3} and edge transitive \cite{Z}, as the hypercube.
However, the balanced hypercube is superior to the hypercube in a sense that it has a smaller diameter than that of the hypercube and supports an
efficient reconfiguration without changing the adjacent relationship among tasks \cite{W}. More desired properties of balanced
hypercubes have been shown in the literature.
Let $BH_{n}$ be an $n$-dimensional balanced hypercube.
Wu and Huang~\cite{H3} proved that $BH_{n}$ is bipartite graph and vertex transitive.
Huang and Wu~\cite{H4} studied resource placement problem in $BH_{n}$.
Zhou et al.~\cite{Z} obtained that $BH_{n}$ is edge transitive and
determined the cyclic connectivity and cyclic edge-connectivity of $BH_{n}$.
Xu et al.~\cite{X} showed that $BH_{n}$ is edge-pancyclic
and Hamiltonian laceable. Yang \cite{Y1} further showed that $BH_{n}$ is
bipanconnected for all $n\geq 1$. Yang~\cite{Y2} determined super connectivity and super edge
connectivity of $BH_{n}$. Yang~\cite{Y3} also gave the conditional diagnosability of $BH_{n}$ under the PMC model.
L\"{u} et al.~\cite{L1} obtained (conditional) matching
preclusion number of $BH_{n}$. L\"{u} and Zhang~\cite{L2} further proved that $BH_{n}$ is hyper-Hamiltonian
laceable. Cheng et al.~\cite{C2} showed that there exist two vertex-disjoint path cover in $BH_{n}$.

Since faulty vertices and faulty edges may happen in the real world networks, the fault-tolerant of networks has been proposed, see \cite{F,H2,T,X0,Y5} etc..
There are some results about the fault-tolerance of $BH_{n}$.
Let $F_{v}$ (resp. $F_{e}$) be the set of faulty vertices (resp. faulty edges) in the $n$-dimensional balanced hypercube $BH_{n}$.
As examples, Zhou et al.~\cite{Z1} derived that $BH_{n}$ is $(2n-2)$-edge-fault Hamiltonian laceable for $n\geq 2$, i.e., if at most $2n-2$ faulty edges occurs, the resulting graph remains Hamiltonian laceable.
Cheng et al. in~\cite{C1} proved that if $|F_{v}|\leq n-1$ faulty vertices occurs, every fault-free edge of $BH_{n}-F_{v}$ lies on a fault-free cycle of every even length
from $4$ to $2^{2n}-2|F_{v}|$.
Cheng et al. in~\cite{C} derived that $BH_{n}$ is $(2n-3)$-edge-fault-bipancyclic, that is, there are at most  $|F_{e}|\leq 2n-3$ faulty edges, every fault-free edge of $BH_{n}-F_{e}$ lies on a fault-free cycle of every even length from $4$ to $2^{2n}$.
Hao et al.~\cite{H} gave the existence of $(2n-2)$-
edge-fault-tolerant cycle embedding in $BH_{n}$.
But up to now, almost all results about fault-tolerance in $BH_{n}$ with only faulty vertices or only faulty edges.

In this paper, we consider cycle embedding of $BH_{n}$ with both faulty vertices and faulty edges,
and prove that there exists a fault-free cycle of length $2^{2n}-2|F_{v}|$ in $BH_{n}$ with $|F_{v}|+|F_{e}|\leq 2n-2$ and $|F_{v}|\leq n-1$ for $n\geq 2$.
Since $BH_{n}$ is a bipartite graph with two
partite sets of equal size, the cycle of length $2^{2n}-2|F_{v}|$ is the longest in the worst-case. This paper improves the results about the longest fault-free cycle in $BH_{n}$ with only faulty vertices or only faulty edges.

The rest of this paper is organized as follows. Section 2 presents some
necessary definitions of graphs and properties of balanced hypercubes as preliminaries.
The proof of our main result is given in Section 3. Section 4 concludes this paper.

\section{Preliminaries}

Let $G=(V,E)$ be a simple undirected graph, where $V$ is the vertex-set of $G$ and $E$ is the edge set of $G$. The {\it neighborhood} of a vertex $u$ in $G$ is the set of all vertices adjacent to $u$ in $V(G)$,
denoted by $N_{G}(u)$. The cardinality $|N_{G}(u)|$ represents the {\it degree} of $u$ in $G$, denoted by $d_{G}(u)$.
A {\it path} $P$, denoted by $P=\langle u_{0},u_{1},\ldots, u_{k}\rangle$, is a sequence of vertices where two successive vertices are adjacent in $G$.
$u_{0}$ and $u_{k}$ are called the {\em end-vertices} of a path $P$. A {\it $u-v$ path} which is denoted by $P[u,v]$ is a path with endpoints $u$ and $v$. For $0\leq i< j\leq k$, let
$P=\langle u_{0},u_{1},\ldots, u_{k}\rangle=\langle u_{0},u_{1},\ldots,u_{i},P[u_{i},u_{j}],u_{j},\ldots,u_{k}\rangle$, where $P[u_{i},u_{j}]=\langle u_{i},u_{i+1},\ldots,u_{j}\rangle$.
A {\it cycle} is a path that two ends are the same vertex. A cycle (resp. path) is a {\it Hamiltonian cycle} (resp. {\it Hamiltonian path}) if it traverses all the vertices of $G$ exactly once.
The {\it length} of a path is the number of edges in it.

A graph $G$ is said to be {\it Hamiltonian connected} if there exists a Hamiltonian path
between any two vertices of $G$. It is easy to see that any bipartite graph with at least three vertices is not
Hamiltonian connected.
A bipartite graph $G$ is {\it Hamiltonian laceable} (resp. {\it strongly Hamiltonian
laceable}) if there exists a of length $|V(G)|-1$ (resp. $|V(G)|-2$) between two arbitrary vertices from different (resp. same) partite sets.
A bipartite graph $G$ with bipartition $V_{0}\cup V_{1}$
is {\it hyper-Hamiltonian laceable} if it is Hamiltonian
laceable and, for any vertex $v \in V_{i}$ ($i \in\{0, 1\})$, there is a Hamiltonian path in $G-v$
between any pair of vertices in $V_{1-i}$.

An {\it $n$-dimensional balanced hypercube}, denoted by $BH_n$, is proposed by Wu and Huang \cite{W} which is defined as an undirected simple graph.

\begin{defi}\label{defi=1}
An $n$-dimensional balanced hypercube $BH_n$ with $n\geq 1$ has $2^{2n}$ vertices with
addresses $(a_0,a_1,\ldots,a_{n-1})$, where $a_i\in \{0,1,2,3\}$ for each $0\leq i\leq n-1$.
Each vertex $(a_0,a_1,\ldots,a_{n-1})$ is adjacent to the following $2n$ vertices:\\
$$\left \{ \begin{array}{l}
((a_0\pm 1)({\rm mod}\ 4), a_1,\ldots, a_{n-1});\\
((a_0\pm 1)({\rm mod}\ 4), a_1,\ldots, a_{j-1},(a_j+(-1)^{a_0})({\rm mod}\ 4),a_{j+1},\ldots,a_{n-1}),
\end{array} \right.$$

where $j$ is an integer with $1\leq j\leq n-1$.
\end{defi}

`$+$' and `$-$' are the operation with $``({\rm mod}\ 4)"$. For notation convenience, $``({\rm mod}\ 4)"$ is omitted.
The balanced hypercube $BH_{n}$ can also be recursively defined.

\begin{defi}\label{defi=2}
$BH_n$ can be recursively constructed as follows:
\begin{enumerate}
\itemsep -1pt
\item [\rm (1)] $BH_1$ is a cycle consisting of $4$ vertices labeled as
$0, 1, 2, 3$ respectively.
\item [\rm (2)] For $n\geq 2$, $BH_{n}$ consists of four copies of $BH_{n-1}$, denoted by $BH_{n-1}^i$ for each $i\in\{0,1,2,3\}$. Each vertex
$(a_0,a_1,\ldots, a_{n-2},i)$ in $BH_{n-1}^i$ has two extra adjacent vertices which are also called extra neighbors:
\begin{enumerate}
\itemsep -1pt
\item [\rm (2.1)] $(a_0\pm1,a_1,\ldots, a_{n-2},i+1)$ in $BH_{n-1}^{i+1}$ if $a_0$ is even.
\item [\rm (2.2)] $(a_0\pm1,a_1,\ldots, a_{n-2},i-1)$ in $BH_{n-1}^{i-1}$ if $a_0$ is odd.
\end{enumerate}
\end{enumerate}

\end{defi}

\begin{figure}[ht]
\begin{center}
\unitlength 4mm
\begin{picture}(20,12)

\put(1, 3){\circle*{0.4}} \put(0.6, 2){$3$}

\put(1,3){\line(1, 0){5}} \put(1,3){\line(0, 1){5}}

\put(6, 3){\circle{0.4}}\put(5.7, 2){$2$}

\put(6,3){\line(0, 1){5}}\put(1,8){\line(1, 0){5}}

\put(1, 8){\circle{0.4}} \put(0.6, 8.3){$0$}

\put(6, 8){\circle*{0.4}}\put(5.7, 8.3){$1$}

\put(10, 1){\circle*{0.4}}\put(10,1){\line(0, 1){9}}\put(10,1){\line(1, 0){3}}\put(10,1){\line(2, 1){6}}
\qbezier(10, 1)(14.5, -1)(19, 1)\put(8, 1){{\footnotesize $(3,3)$}}\put(8, 4){{\footnotesize $(0,3)$}}
\put(8, 7){{\footnotesize $(3,0)$}}\put(8, 10){{\footnotesize $(0,0)$}}

\put(10, 4){\circle{0.4}}\put(10,4){\line(1, 0){9}}\put(10,4){\line(1, 2){3}}

\put(10, 7){\circle*{0.4}}\put(10,7){\line(1, 0){9}}\put(10,7){\line(1, -2){3}}

\put(10, 10){\circle{0.4}}\put(10,10){\line(2, -1){6}}\put(10,10){\line(1, 0){3}}
\qbezier(10, 10)(14.5, 12)(19, 10)

\put(13, 1){\circle{0.4}} \put(13,1){\line(0, 1){3}}\qbezier(13, 1)(14.5, 5.5)(13, 10)
\put(12.5, 0){{\footnotesize $(2,3)$}}
\put(15, 0){{\footnotesize $(3,2)$}}

\put(13, 4){\circle*{0.4}}\put(13,4){\line(2, -1){6}}\put(12.5, 4.2){{\footnotesize $(1,3)$}}
\put(15, 4.2){{\footnotesize $(0,2)$}}

\put(13, 7){\circle{0.4}}\put(13,7){\line(0, 1){3}} \put(13,7){\line(2, 1){6}}\put(12.5, 6){{\footnotesize $(2,0)$}}
\put(15, 6){{\footnotesize $(3,1)$}}

\put(13, 10){\circle*{0.4}} \put(12.5, 10.5){{\footnotesize $(1,0)$}}\put(15, 10.5){{\footnotesize $(0,1)$}}

\put(16, 1){\circle*{0.4}}\put(16,1){\line(0, 1){3}}\qbezier(16, 1)(14.5, 5.5)(16, 10)\put(16,1){\line(1, 2){3}}

\put(16, 4){\circle{0.4}}

\put(16, 7){\circle*{0.4}}\put(16,7){\line(0, 1){3}}

\put(16, 10){\circle{0.4}}\put(16,10){\line(1, 0){3}} \put(16,10){\line(1, -2){3}}

\put(19, 1){\circle{0.4}}\put(19,1){\line(0, 1){9}}\put(19,1){\line(-1, 0){3}}

\put(19, 4){\circle*{0.4}}

\put(19, 7){\circle{0.4}}

\put(19, 10){\circle*{0.4}}\put(19.2, 1){{\footnotesize $(2,2)$}}\put(19.2, 4){{\footnotesize $(1,2)$}}
\put(19.2, 7){{\footnotesize $(2,1)$}}\put(19.2, 10){{\footnotesize $(1,1)$}}
\end{picture}
\end{center}\vspace{-0.3cm}
\caption{$BH_1$ and $BH_2$} \label{F-1}
\end{figure}
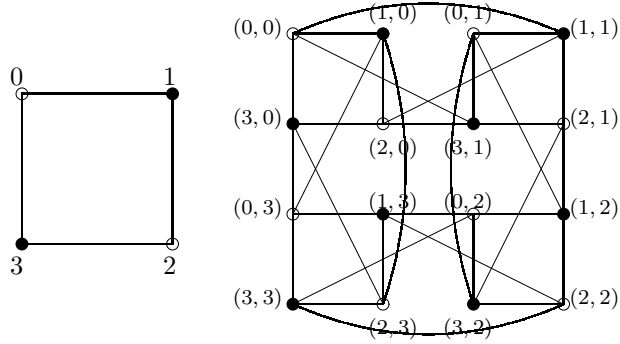

$BH_1$ and $BH_2$ are illustrated in Figure~\ref{F-1}.

The first element $a_{0}$ of vertex $(a_0,a_1,\ldots, a_{n-2},a_{n-1})$ is called the {\it inner index}, and the other elements $a_{i}$,
for all $1\leq i \leq n-1$, are called {\it $i$-dimensional index}. By~\cite[Lemma 5]{Y4},
for each vertex $u$ of $BH_n$, there exists a unique vertex, denoted by $u^{\ast}$, such that $N_{BH_n}(u)=N_{BH_n}(u^{\ast})$.

Since $BH_{n}$ is a bipartite graph,
we refer to a vertex with an odd inner index as a {\it black vertex} and a vertex with an even inner index as a {\it white vertex}.
If two adjacent vertices $u$ and $v$ differ in only the inner index, the edge $e=(u,v)$ is said to be {\it $0$-dimensional} and $v$ is a {\it $0$-dimensional neighbor} of $u$. Let $[n]=\{0,1,\ldots,n-1\}$ (but not general meaning $\{1,2,\ldots,n\}$).
Likewise if two adjacent vertices $u$ and $v$ not only differ in the inner index,
but also differ in some $j$-dimensional index ($j\in [n]$), the edge $e=(u,v)$ is said to be {\it $j$-dimensional}, and $v$ is the {\it $j$-dimensional
neighbor} of $u$.
Let $E_{j}$, where $j\in [n]$, denote the set of all $j$-dimensional edges. For $i\in \{0,1,2,3\}$ and $j\in [n]$, we use $BH_{n-1}^{j,i}$
to denote $(n-1)$-dimensional sub-balanced hypercubes of the $BH_{n}$ induced by all vertices labeled by $(a_0,a_1,\ldots,a_{j-1},i,a_{j+1},\ldots, a_{n-2},a_{n-1})$.
Obviously, $BH_{n}-E_{j}=\bigcup\limits_{i=0}^{3}BH_{n-1}^{j,i}$ and $BH_{n-1}^{j,i}\cong BH_{n-1}$.
The edges between $BH_{n-1}^{j,i}$ and $BH_{n-1}^{j,i+1}$ are called {\it $j$-dimensional crossing edges}, where $i\in \{0,1,2,3\}$ and $j\in \{1,2,\ldots,n-1\}$.
If $j=n-1$, $BH_{n-1}^{j,i}$ and $E_{j}$ are denoted by $BH_{n-1}^{i}$ and $E_{c}$, respectively, where $i\in \{0,1,2,3\}$
and $(n-1)$-dimension crossing edges are called {\it crossing edges}.
Let $E_{0}$ be the set of all $0$-dimensional edges,
L\"{u} et al. \cite{L1} obtained that $BH_{n}-E_{0}$ has four components, each component is isomorphic to $BH_{n-1}$.

Some known basic properties of $BH_n$ are given as follows.

\begin{lem}{\rm(\cite{W})}\label{vertex-trans}
 The balanced hypercube $BH_{n}$ is bipartite and vertex transitive.
\end{lem}

\begin{lem}\label{edge-trans}
{(\rm\cite{Z})} The balanced hypercube $BH_{n}$ is edge transitive.
\end{lem}

\begin{lem}{\rm(\cite{X})}\label{H-L}
The balanced hypercube $BH_{n}$ is edge-bipancyclic and Hamiltonian laceable for all $n\geq 1$.
\end{lem}

\begin{lem}{(\rm\cite{Y1})}\label{bipan}
The balanced hypercube $BH_{n}$ is bipanconnected for all $n\geq 1$.
\end{lem}

\begin{lem}{\rm(\cite{L2})}\label{laceable}
The balanced hypercube $BH_{n}$ is hyper-Hamiltonian laceable for $n\geq1$. Thus, $BH_{n}$ is strongly Hamiltonian laceability.
\end{lem}

\begin{lem}{\rm(\cite{C2})}\label{two-node}
Let $X$ and $Y$ be two distinct partite sets of $BH_{n}$. Assume $u$
and $x$ are two different vertices in $X$, $v$ and $y$ are two different vertices in $Y$.
Then there exist two vertex-disjoint paths $P[x, y]$ and $R[u, v]$, and $V(P[x,y])\cup V(R[u,v])= V(BH_{n})$.
\end{lem}

\begin{lem}\label{8-cycle}
\begin{enumerate}
\item [{\rm (1)}]{\rm(\cite{X})}
Let $e=(u,v)$ be an edge of $BH_{n}$. $e$ is contained in a cycle $C$ of length $8$ in $BH_{n}$ such that $|E(C)\cap E(BH_{n-1}^{i})|=1$, where $i=0,1,2,3$.
\item [{\rm (2)}]\ {\rm(\cite{C1})}
Let $e=(u,v)$ be an edge of $BH_{n}$. If $e$ is along dimension $j$ ($1\leq j\leq n-1$), then there are
$2n-2$ internal vertex-disjoint paths of length $7$ joining $u$ and $v$ such that each path has only one edge in each $BH_{n-1}^{j,i}$, where $i\in\{0,1,2,3\}$.

\item [{\rm (3)}]\ {\rm(\cite{C1})} Let $e=(x,y)$ be any edge in $BH_{n-1}^{j,0}$, then there exist two internal vertex-disjoint paths
$\langle x,y_{1},x_{1},y_{2},x_{2},y_{3},x_{3},y\rangle$ and $\langle x,y_{1}',x_{1}',y_{2}',x_{2}',\\y_{3}',x_{3}',y\rangle$ in $BH_{n}$ such that
$(x_{i},y_{i}),(x_{i}',y_{i}')\in E(BH_{n-1}^{j,i})$, where $1\leq j\leq n-1$ and $i=1,2,3$.
\end{enumerate}
\end{lem}

\begin{lem}{\rm(\cite{H})}\label{edge-dis}
Let $A=\{e_{1},e_{2},\ldots,e_{k}\}$ for $k\geq2$ be the set of $j$-dimensional edges in $BH_{n-1}^{i}$, where $j\in[n]$ and $i\in \{0,1,2,3\}$.
Then there are cycles $D_{1},D_{2},\ldots,D_{k}$ which contain $e_{1},e_{2},\ldots,e_{k}$, respectively, such that the length of each cycle
is $8$, $D_{1},D_{2},\ldots,D_{k}$ are edge disjoint and $|E(D_{m})\cap E(BH_{n-1}^{t})|=1$ for every $m\in \{1,2,\ldots,k\}$ and $t\in \{0,1,2,3\}$.
Moreover, for any $i,j\in[k]$, if $e_{i}\cap e_{j}=\emptyset$, then  $D_{i},D_{j}$ are also vertex-disjoint.
If $e_{i}\cap e_{j}=\{u\}$, then $D_{i},D_{j}$ are vertex-disjoint except a vertex $u$.

\end{lem}

For the fault tolerance of $BH_{n}$, Cheng et al.~\cite{C} and~\cite{C1} derived the following two results, respectively..

\begin{lem} {\rm(\cite{C})}\label{edge-2n-3}
Let $F_{e}$ be the set of faulty vertices in $BH_{n}$ with $|F_{e}|\leq 2n-3$. Then every fault-free edge of $BH_{n}-F_{e}$ lies on a fault-free
cycle of every even length from $4$ to $2^{2n}$ inclusive, where $n\geq2$.
\end{lem}

\begin{lem} {\rm(\cite{C1})}\label{vertex-n-1}
Let $F_{v}$ be the set of faulty vertices in $BH_{n}$ with $|F_{v}|\leq n-1$. Then every fault-free edge of $BH_{n}-F_{v}$ lies on a fault-free
cycle of every even length from $4$ to $2^{2n}-2|F_{v}|$ inclusive, where $n\geq1$.
\end{lem}

Hao et al.~\cite{H} and Zhou et al.~\cite{Z1} obtained the following two results, respectively.

\begin{lem}{\rm(\cite{H})}\label{H-2n-2}
Let $F_{e}$ be a faulty edge set of the balanced hypercube $BH_{n}$ with $|F_{e}|\leq 2n-2$, $n\geq 1$, then there exists a fault-free
Hamiltonian path between any two adjacent vertices $x$ and $y$ in $BH_{n}$.
\end{lem}

\begin{lem}{\rm(\cite{Z1})}\label{Z}
$BH_{n}$ is $(2n-2)$-edge-fault Hamiltonian laceable for $n\geq 2$.
\end{lem}

\section{Main results}

In this section, we consider cycle embedding in faulty balanced hypercube $BH_{n}$. For this purpose, we need the following lemma.

\begin{lem}\label{n-1-laceable}
Let $F_{v}$ (resp. $F_{e}$) be the set of faulty vertices (resp. faulty edges) with $|F_{v}|+|F_{e}|\leq n-1$ in $BH_{n}$ for $n\geq 1$. Then,
there is a fault-free path in $BH_{n}$ whose length is $2^{2n}-2|F_{v}|-1$ between any two adjacent vertices $x$ and $y$ in $BH_{n}$.
\end{lem}

\f {\bf Proof.} Assume $x$ is a white vertex and $y$ is a black vertex and $(x,y)\in E(BH_{n})$. Recall that $f_{e}=|F_{e}|$, $f_{v}=|F_{v}|$, $f=f_{e}+f_{v}$.  We only need to prove the result holds for $f=n-1$.
 We show the lemma by induction on $n$. If $n=1$, $f=0$, $BH_{1}$ is a cycle of length $4$, so the result is right. Now we consider $n=2$. If $f_{e}=0$, then $f_{v}=1$.
Since $x$ and $y$ are adjacent, regard $(x,y)$ as a fault-free edge, by Lemma~\ref{vertex-n-1}, the edge $(x,y)$ lies in a fault-free cycle of length $2^{4}-2$,
this implies there exists a fault-free path $P[x,y]$ of length $2^{4}-2-1=13$.
If $f_{v}=0$, then $f_{e}=1$.
Since $x$ and $y$ are adjacent, by Lemma~\ref{H-2n-2}, there exists a fault-free Hamiltonian $x-y$ path of length $2^{4}-1=15$.
The lemma holds for $n=2$.

Assume that the lemma holds for $BH_{n-1}$ for $n\geq 3$. Now we consider $BH_{n}$.
If $BH_{n}$ has no faulty edges, then $f_{v}\leq n-1$. By Lemma~\ref{vertex-n-1},
there exists a fault-free path of length $2^{2n}-2f_{v}-1$ between $x$ and $y$.
So we assume that there is at least one faulty edge. Without loss of generality, suppose there exists an $(n-1)$-dimensional faulty edge.
We can divide $BH_{n}$ into $BH_{n-1}^{i}$, $i=0,1,2,3$, along $(n-1)$-dimension.
For $i\in \{0,1,2,3\}$, let $f_{e}^{i}=|F_{e}\cap E(BH_{n-1}^{i})|$ and $f_{v}^{i}=|F_{v}\cap V(BH_{n-1}^{i})|$.
Moreover, we have that $0\leq f_{e}^{i}+f_{v}^{i}\leq n-2$ for each $i\in \{0,1,2,3\}$. There are the following two cases.

Case 1. $x$ and $y$ are in the same $BH_{n-1}^{i}$ for some $i\in\{0,1,2,3\}$. Without loss of generality, let $i=0$.

By the inductive hypothesis, there exists a fault-free $x-y$ path $P_{0}$ in $BH_{n-1}^{0}$ of length $2^{2(n-1)}-2f_{v}^{0}-1$.
We claim that there exists at least one $k\in \{0,1,\ldots,n-2\}$ such that $P_{0}$ contains at least $2n-1$ $k$-dimensional
edges. Otherwise, for every $k\in \{0,1,\ldots,n-2\}$, $P_{0}$ contains at most $2n-2$ $k$-dimensional edges. As a result, $P_{0}$ contains at most $(n-1)(2n-2)$ edges.
Since $|P_{0}|=2^{2(n-1)}-2f_{v}^{0}-1\geq 2^{2(n-1)}-2n+1> (n-1)(2n-2)$ for $n\geq 3$, it is a contradiction.
Since $f=n-1\leq 2n-1$ for $n\geq 3$, by Lemma~\ref{edge-dis}, we can find a $k$-dimensional edge $(u_{0},v_{0})$ (assume $u_{0}$ is white, $v_{0}$ is black) on $P_{0}$,
such that $(u_{0},v_{0})$ is contained in a fault-free $8$-cycle $C$, say $C=\langle u_{0},v_{0},u_{3},v_{3},u_{2},v_{2},u_{1},v_{1},u_{0}\rangle$,
where $(u_{i},v_{i})\in E(BH_{n-1}^{i})$ for $i\in \{0,1,2,3\}$.

By the inductive hypothesis in $BH_{n-1}^{i}$, for each $i\in\{1,2,3\}$, there exists a fault-free $u_{i}-v_{i}$ path $P_{i}$ in $BH_{n-1}^{i}$ of length $2^{2(n-1)}-2f_{v}^{i}-1$.
Therefore, $\langle x,P[x,u_{0}],u_{0},v_{1},P_{1},u_{1},v_{2},P_{2},\\u_{2},v_{3},P_{3},u_{3},v_{0},P[v_{0},y],y\rangle$ is a desired fault-free $x-y$ path in $BH_{n}$ of length
$\Sigma_{i=0}^{3}(2^{2(n-1)}-2f_{v}^{i}-1)+3=2^{2n}-2f_{v}-1$.

Case 2. $x$ and $y$ are in two distinct $BH_{n-1}^{i}$'s, where $i\in \{0,1,2,3\}=[4]$.

Without loss of generality, suppose $x\in V(BH_{n-1}^{0})$. Since $x$ and $y$ are adjacent, $y\in V(BH_{n-1}^{1})$.
By Lemma~\ref{8-cycle}(2) and $f\leq n-1\leq 2n-2$ for $n\geq 3$, there exists a fault-free $x-y$ path $P$, say $P=\langle x,v_{0},u_{3},v_{3},u_{2},v_{2},u_{1},y\rangle$, of length $7$
such that $P$ has only one edge in $BH_{n-1}^{i}$ for each $i=0,1,2,3$.
By the inductive hypothesis, there exists a fault-free $x-v_{0}$ (resp. $u_{3}-v_{3}$, $u_{2}-v_{2}$ and $u_{1}-y$) path $P_{0}$ (resp. $P_{3}$, $P_{2}$ and $P_{1}$) in $BH_{n-1}^{0}$
(resp. $BH_{n-1}^{3}$, $BH_{n-1}^{2}$ and $BH_{n-1}^{1}$) of length $2^{2(n-1)}-2f_{v}^{0}-1$ (resp. $2^{2(n-1)}-2f_{v}^{3}-1$, $2^{2(n-1)}-2f_{v}^{2}-1$ and $2^{2(n-1)}-2f_{v}^{1}-1$).
A fault-free $x-y$ path of $BH_{n}$ can be constructed as $\langle x,P_{0},v_{0},u_{3},P_{3},v_{3},u_{2},P_{2},v_{2},u_{1},P_{1},y\rangle$ of length
$\Sigma_{i=0}^{3}(2^{2(n-1)}-2f_{v}^{i}-1)+3=2^{2n}-2f_{v}-1$.

\medskip
By the above cases, the proof is complete.
\hfill\qed

\begin{lem}\label{n=2}
Let $v$ be a faulty vertex and $e$ be a faulty edge in $BH_{2}$.
Then there exists a fault-free cycle of length $14$ in $BH_{2}$.
\end{lem}

\f {\bf Proof.} By Lemma~\ref{edge-trans}, $BH_{2}$ is edge transitive, we may assume that $a=((0,0),(1,1))$ is a faulty edge.
For each possible faulty vertex $v\in V(BH_{2})$, we can find a desired cycle. All the cases are listed in Table~\ref{table=1}.
\hfill\qed

\begin{table}[!htbp]
\begin{center}
\begin{tabular}{>{\scriptsize}l}
\hline
$v=(0,0)$ or $(1,0)$ \\
$\langle (3,0),(2,0),(3,1),(0,1),(1,1),(2,1),(1,2),(2,2),(3,2),(0,2),(1,3),(2,3),(3,3),(0,3),(3,0)\rangle$\\

$v=(3,0)$ or $(2,0)$ \\
$\langle (0,0),(1,0),(0,3),(3,3),(2,3),(1,3),(0,2),(3,2),(2,2),(1,2),(2,1),(1,1),(0,1),(3,1),(0,0)\rangle$\\

$v=(0,1)$ or $(1,1)$ \\
$\langle (0,0),(1,0),(2,0),(3,1),(2,1),(1,2),(2,2),(3,2),(0,2),(1,3),(2,3),(3,3),(0,3),(3,0),(0,0)\rangle$\\

$v=(3,1)$ or $(2,1)$ \\
$\langle (0,0),(1,0),(2,0),(1,1),(0,1),(3,2),(2,2),(1,2),(0,2),(1,3),(2,3),(3,3),(0,3),(3,0),(0,0)\rangle$\\

$v=(0,3)$ or $(1,3)$ \\
$\langle (0,0),(3,0),(2,0),(1,0),(2,3),(3,3),(2,2),(3,2),(0,2),(1,2),(2,1),(1,1),(0,1),(3,1),(0,0)\rangle$\\

$v=(3,3)$ or $(2,3)$ \\
$\langle (0,0),(3,0),(2,0),(1,0),(0,3),(1,3),(0,2),(3,2),(2,2),(1,2),(2,1),(1,1),(0,1),(3,1),(0,0)\rangle$\\

$v=(0,2)$ or $(1,2)$ \\
$\langle (0,0),(1,0),(2,0),(3,0),(0,3),(3,3),(2,3),(1,3),(2,2),(3,2),(0,1),(1,1),(2,1),(3,1),(0,0)\rangle$\\

$v=(3,2)$ or $(2,2)$ \\
$\langle (0,0),(1,0),(2,0),(3,0),(0,3),(3,3),(2,3),(1,3),(0,2),(1,2),(2,1),(1,1),(0,1),(3,1),(0,0)\rangle$\\

\hline
\end{tabular}

\vskip 0.2cm
\caption{{\small  Desired cycles of length $14$ in $BH_{2}$ with a faulty vertex $v$ and a faulty edge ((0,0),(1,1))}}\label{table=1}
\end{center}
\end{table}

The following is our main result.

\begin{theorem}\label{main}
Let $F_{v}$ and $F_{e}$ be the set of faulty vertices and faulty edges, respectively,
 with $|F_{v}|\leq n-1$ and $|F_{v}|+|F_{e}|\leq 2n-2$ in $BH_{n}$, where $n\geq 2$.
Then there exists a fault-free cycle of length $2^{2n}-2|F_{v}|$ in $BH_{n}$.
\end{theorem}

\f {\bf Proof.} We prove this theorem by induction on $n$. Recall that $E_{j}$ is the set of $j$-dimensional edges in $BH_{n}$ for $j\in[n]$. Since $|F_{v}|+|F_{e}|\leq 2n-2$, there exists a $j\in [n]$ such that $|F_{e}\cap E_{j}|\leq1$.
Without loss of generality, let $j=n-1$. Then we can divide $BH_{n}$ into four components, say $BH_{n-1}^{i}$, $i=0,1,2,3$, along $(n-1)$-dimension.
Let $F_{c}=F_{e}\cap E_{n-1}$, then $f_{c}=|F_{c}|\leq 1$.
Recall that $f_{v}=|F_{v}|$, $f_{e}=|F_{e}|$, $f=f_{v}+f_{e}$, $f_{e}^{i}=|F_{e}\cap E(BH_{n-1}^{i})|$ and $f_{v}^{i}=|F_{v}\cap V(BH_{n-1}^{i})|$, so $f_{i}=f_{e}^{i}+f_{v}^{i}$ for $i\in \{0,1,2,3\}$.

If $n=2$, then $f_{v}+f_{e}\leq 2$ and $f_{v}\leq 1$.
If $f_{v}=1$ and $f_{e}=0$, then by Lemma~\ref{vertex-n-1}, the result is right.
If $f_{v}=1$ and $f_{e}=1$, the result follows from Lemma~\ref{n=2}.
If $f_{v}=0$ and $f_{e}\leq 2$. Choose a fault-free edge $(u,v)$, by Lemma~\ref{H-2n-2}, there
exists a fault-free Hamiltonian path $P$ between $u$ and $v$ in $BH_{2}$, then $\langle u,P,v,u\rangle$ is the desired cycle of $BH_{2}$.
The theorem is true for $n=2$.

For $n\geq3$, we assume that the theorem is true in $BH_{m}$ with $f_{v}+f_{e}\leq 2m-2$ and $f_{v}\leq m-1$ for every integer $2\leq m< n$. Now we consider $BH_{n}$ as follows.
If $BH_{n}$ has no faulty vertices, then $f_{v}=0$ and $f_{e}\leq2n-2$. Choose a fault-free edge $(x,y)$, since $f_{e}\leq 2n-2$, by Lemma~\ref{H-2n-2},
there exists a fault-free Hamiltonian path $P$ from $x$ to $y$ in $BH_{n}$. Let $C=\langle x,P,y,x\rangle$. Clearly, $|E(C)|=2^{2n}=2^{2n}-2f_{v}$. So we assume $1\leq f_{v}\leq n-1$ in the following. We only need to consider the following three cases.

\medskip

Case 1.  For each $i\in \{0,1,2,3\}$, $0\leq f_{i}\leq2n-4$.

Subcase 1.1. For each $i\in \{0,1,2,3\}$, $0\leq f_{i}\leq n-2$.

By the inductive hypothesis, there exists a cycle $C_{0}$ in $BH_{n-1}^{0}$ of length $2^{2(n-1)}-2f_{v}^{0}$. Select an $(u_{0},v_{0})\in E(C_{0})$ and assume $u_{0}$ is white.
By Lemma~\ref{8-cycle}(3), since $f_{c}\leq 1$, $(u_{0},v_{0})$ is contained in two $8$-cycles which has only one common edge $(u_{0},v_{0})$, so we can find a
$8$-cycle $C$, say $C=\langle u_{0},v_{0},u_{3},v_{3},u_{2},v_{2},u_{1},v_{1},u_{0}\rangle$,
where $(u_{i},v_{i})\in E(BH_{n-1}^{i})$ for $i\in \{0,1,2,3\}$, such that the crossing edge $(u_{i},v_{i+1})$ for each $i\in \{0,1,2,3\}$ is fault-free, where $v_{4}=v_{0}$.
Since $f_{i}\leq n-2$ and $BH_{n-1}^{i}\cong BH_{n-1}$ for $i\in\{1,2,3\}$, by Lemma~\ref{n-1-laceable}, there exists a fault-free $u_{i}-v_{i}$ path $P_{i}$ of length $2^{2(n-1)}-2f_{v}^{i}-1$ in $BH_{n-1}^{i}$.
Let $P_{0}=C_{0}-(u_{0},v_{0})$. So $\langle u_{0},v_{1},P_{1},u_{1},v_{2},P_{2},u_{2},v_{3},P_{3},u_{3},v_{0},\\P_{0},u_{0}\rangle$ is a desired fault-free cycle in $BH_{n}$ of length
$2^{2n}-2f_{v}$.

Subcase 1.2. There exists only one $i\in \{0,1,2,3\}$ such that $n-1\leq f_{i}\leq 2n-4$. Without loss of generality, let $i=0$.
It implies that $0\leq f_{j}\leq n-2$ for all $j\in\{1,2,3\}$.

Subcase 1.2.1. $n-1\leq f_{0}\leq 2n-4$, $f_{v}^{0}\leq n-2$.

By the inductive hypothesis, there exists a fault-free cycle $C_{0}$ in $BH_{n-1}^{0}$ of length $2^{2(n-1)}-2f_{v}^{0}$.
By the similar discussions as Case 1 in the proof of Lemma~\ref{n-1-laceable}, we can find an
edge $(u_{0},v_{0})$ (assume $u_{0}$ is white, $v_{0}$ is black) on $C_{0}$,
such that $(u_{0},v_{0})$ is contained in a fault-free $8$-cycle $C$, say $C=\langle u_{0},v_{0},u_{3},v_{3},u_{2},v_{2},u_{1},v_{1},u_{0}\rangle$,
where $(u_{i},v_{i})\in E(BH_{n-1}^{i})$ for $i\in \{1,2,3\}$. Since $0\leq f_{i}\leq n-2$ for $i\in\{1,2,3\}$, by Lemma~\ref{n-1-laceable} in $BH_{n-1}^{i}$,
 there exists a fault-free $u_{i}-v_{i}$ path $P_{i}$ in $BH_{n-1}^{i}$ of length $2^{2(n-1)}-2f_{v}^{i}-1$. Let $P_{0}=C_{0}-(u_{0},v_{0})$.
Therefore, $\langle u_{0},v_{1},P_{1},u_{1},v_{2},P_{2},u_{2},v_{3},P_{3},u_{3},v_{0},P_{0},u_{0}\rangle$ is a desired fault-free cycle in $BH_{n}$ of length
$2^{2n}-2f_{v}$.

Subcase 1.2.2. $n-1\leq f_{0}\leq 2n-4$, $f_{v}^{0}=n-1$.

It implies that $f_{v}^{j}=0$ and $f_{e}^{j}\leq f-f_{v}\leq n-1$ for each $j\in \{1,2,3\}$. We regard one faulty vertex in $BH_{n-1}^{0}$ as a fault-free vertex, so $f_{v}'=f_{v}^{0}-1=f_{v}-1=n-2$.
By the inductive hypothesis, there exists a cycle $C_{0}$ in $BH_{n-1}^{0}$ of length $2^{2(n-1)}-2f_{v}'$, where $f_{v}'=f_{v}^{0}-1$. Note that $C_{0}$ contains at most
one faulty vertex.

If $C_{0}$ does not contain any faulty vertex, choose an edge, say $(u_{0},v_{0})\in E(C_{0})$
(assume $u_{0}$ is white, $v_{0}$ is black).
By Lemma~\ref{8-cycle}(3) and $f_{c}\leq 1$,
$(u_{0},v_{0})$ is contained in a cycle $C$, say $C=\langle v_{0},u_{0},v_{1},u_{1},v_{2},u_{2},v_{3},u_{3},v_{0}\rangle$, of length $8$
in $BH_{n}$ such that $E(C)\cap E(BH_{n-1}^{i})=\{(u_{i},v_{i})\}$ and $(u_{i},v_{i+1})$ is fault-free for each $i\in \{0,1,2,3\}$.
Since for any $i\in \{1,2\}$, $f_{v}^{i}=0$ and $f_{e}^{i}\leq n-1\leq 2n-2$ for $n\geq 3$, by Lemma~\ref{Z},
there exists a fault-free path $P_{i}$ between $u_{i}$ and $v_{i}$ in $BH_{n-1}^{i}$ whose length is $2^{2(n-1)}-1$.
Since $f_{v}^{3}=0$ and $f_{e}^{3}\leq n-1\leq 2n-3$ for $n\geq 3$, by Lemma~\ref{edge-2n-3}, there exists a fault-free cycle $C_{3}$ of length $2^{2(n-1)}-2$ passing through $(u_{3},v_{3})$.
Let $P_{i}=C_{i}-(u_{i},v_{i})$ for $i=0,3$.
A desired fault-free cycle of $BH_{n}$ can be constructed
as $\langle u_{0},v_{1},P_{1},u_{1},v_{2},P_{2},u_{2},v_{3},P_{3},u_{3},v_{0},P_{0},u_{0}\rangle$
whose length is $2^{2(n-1)}-2(f_{v}^{0}-1)-1+2(2^{2(n-1)}-1)+ (2^{2(n-1)}-3)+4=2^{2n}-2f_{v}^{0}=2^{2n}-2f_{v}$.

\begin{figure}

\begin{minipage}[t]{0.5\linewidth}
\centering
\includegraphics[height=5cm,width=5cm]{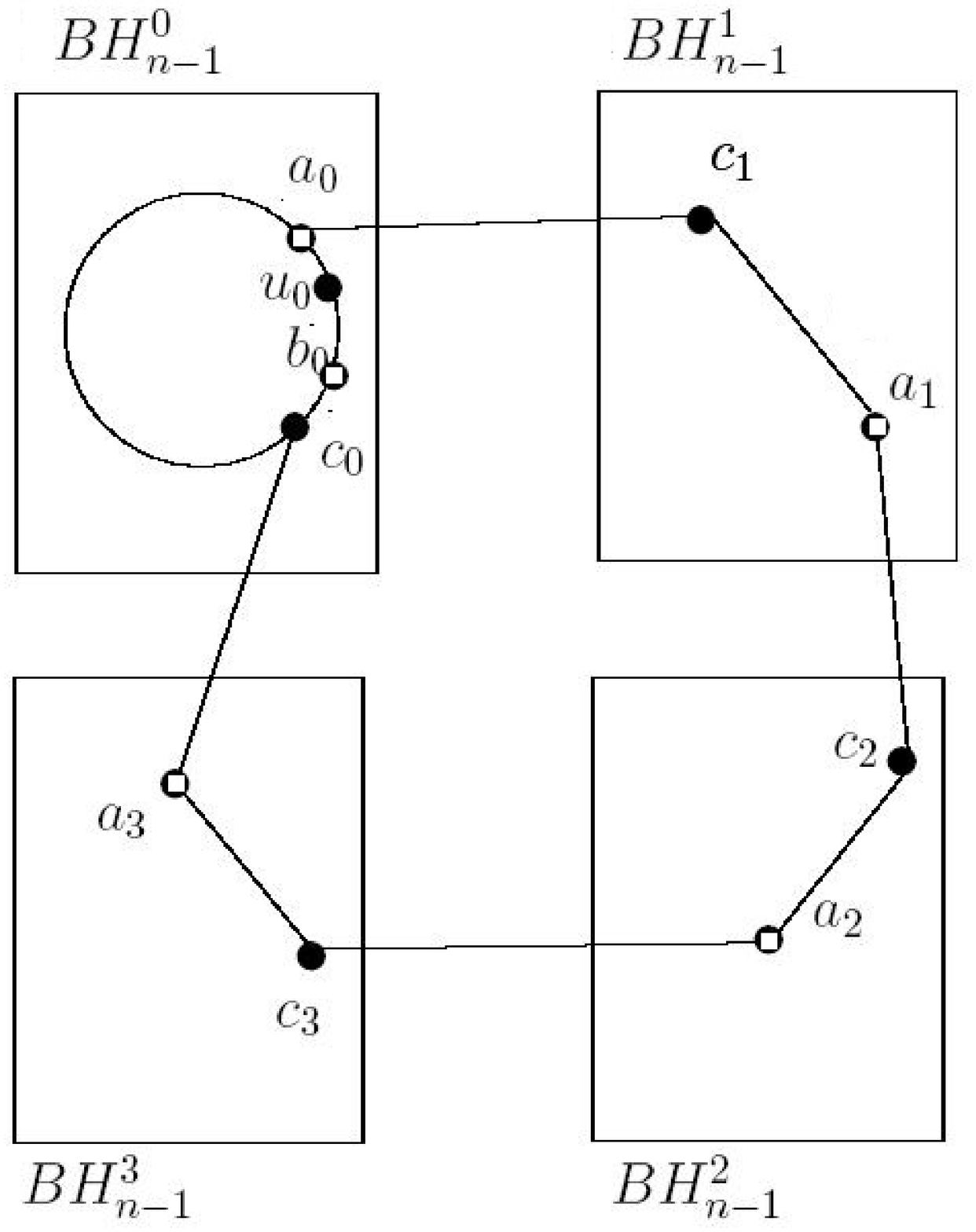}
\vskip-0.2cm
\caption{Illustration of 1.2.2}
\label{1.2}
\end{minipage}
\begin{minipage}[t]{0.5\linewidth}
\centering
\includegraphics[height=5cm,width=5cm]{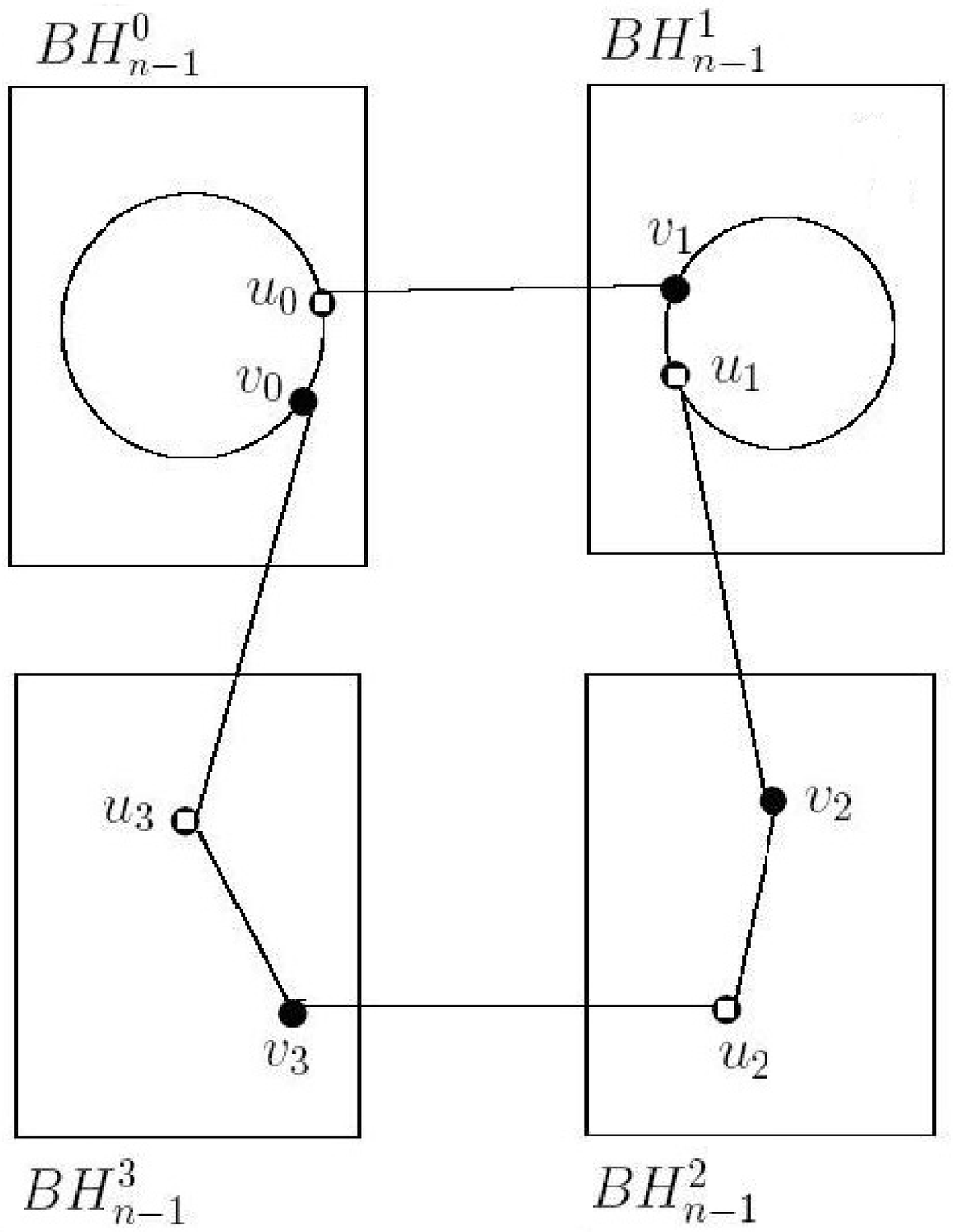}
\vskip-0.6cm
\hskip-1cm\caption{Illustration of 1.3.1.1}
\label{1.3.1.1}
\end{minipage}

\end{figure}

Now assume that $C_{0}$ contains one faulty vertex $u_{0}$ and $u_{0}$ is black. (If $C_{0}$ contains a white faulty vertex,
the discussions are similar, so it is omitted.)
Let $a_{0}$ and $b_{0}$ be two vertices adjacent to $u_{0}$ on $C_{0}$. Let $c_{0}\in C_{0}$ be a neighbor of $b_{0}$.
Recall that each vertex has two extra neighbors and $f_{c}\leq 1$, we can choose the extra neighbors $c_{1}\in V(BH_{n-1}^{1})$ and $a_{3}\in V(BH_{n-1}^{3})$ of $a_{0}$ and $c_{0}$, respectively,
such that both $(a_{0},c_{1})$ and $(c_{0},a_{3})$ are fault-free. For $i\in\{1,2\}$, let $(a_{i},c_{i+1})\in E(BH_{n})$ be a fault-free edge, where $a_{i}\in V(BH_{n-1}^{i})$ be a white vertex
and $c_{i+1}\in V(BH_{n-1}^{i+1})$ be a black vertex. For $i\in \{1,2,3\}$, $f_{e}^{i}\leq n-1\leq 2n-4$ for $n\geq 3$ and $f_{v}^{i}=0$, by Lemma~\ref{Z}, there exists a fault-free path $P_{i}$ of length
$2^{2(n-1)}-1$ between $a_{i}$ and $c_{i}$ in $BH_{n-1}^{i}$. Let $P_{0}=C_{0}-\{u_{0},b_{0}\}$. A desired fault-free cycle of $BH_{n}$ can be constructed as
$\langle a_{0},P_{0},c_{0},a_{3},P_{3},c_{3},a_{2},P_{2},c_{2},a_{1},P_{1},c_{1},a_{0}\rangle$
of length $2^{2n}-2f_{v}$ (see Figure~\ref{1.2}).

Subcase 1.3. There exist two distinct $i,k\in \{0,1,2,3\}$ such that $n-1\leq f_{i}\leq 2n-4$ and $n-1\leq f_{k}\leq 2n-4$. Without loss of generality, assume that $i=0$. Since $f\leq 2n-2$, it implies that $f_{0}=f_{k}=n-1$, $f_{c}=|F_{c}|=0$ and $f_{t}=0$ for any $t\in \{1,2,3\}-\{k\}$.

Subcase 1.3.1. $f_{v}^{0}\leq n-2$ and $f_{v}^{k}\leq n-2$.

By the inductive hypothesis, there exists a cycle $C_{0}$ (resp. $C_{k}$) in $BH_{n-1}^{0}$ (resp. $BH_{n-1}^{k}$) of length $2^{2(n-1)}-2f_{v}^{0}$ (resp. $2^{2(n-1)}-2f_{v}^{k}$).

Subcase 1.3.1.1. $k=1$ or $k=3$.

Since the discussion for the two situations are similar, we may assume $k=1$.

Note that $C_{0}$ has $2^{2(n-1)}-2f_{v}^{0}$ vertices, half is black and half is white, so there are $2^{2n-3}-f_{v}^{0}$ white vertices.
We claim that there exist at least one white vertex in $C_{0}$ such that its extra neighbor in $C_{1}$.
Since for each white vertex, say $u$, in $C_{0}$, there exists at most one other vertex $u^{\ast}$ in $C_{0}$, which has the same extra neighbors as $u$.
So vertices in $C_{0}$ have at least $(2^{2n-3}-f_{v}^{0})/2$ extra neighbors in $C_{1}$. By $(2^{2n-3}-f_{v}^{0})/2>n-1=f_{1}$ because of $f_{v}^{0}\leq n-2$ and $n\geq 3$, there exists a white vertex $u_{0}\in V(C_{0})$ such that  $u_{0}$ has an extra neighbor, say $v_{1}$, on $C_{1}$.
See Figure~\ref{1.3.1.1}.

Let $v_{0}$ be a neighbor of $u_{0}$ in $C_{0}$ and $u_{1}$ be a neighbor of $v_{1}$ in $C_{1}$.
Let $u_{3}\in V(BH_{n-1}^{3})$ (resp. $v_{2}\in V(BH_{n-1}^{2})$) be an extra neighbor of $v_{0}$ (resp. $u_{1}$).
Let $(v_{2},u_{3})\in E(BH_{n})$,
where $v_{2}\in V(BH_{n-1}^{2})$ be a black vertex and $u_{3}\in V(BH_{n-1}^{3})$ be a white vertex.
Since $f_{t}=0$ for $t\in\{2,3\}$, by Lemma~\ref{H-L}, there exists a fault-free path $P_{t}$ of length
$2^{2(n-1)}-1$ between $u_{t}$ and $v_{t}$ in $BH_{n-1}^{t}$. Let $P_{i}=C_{i}-(u_{i},v_{i})$ for $i=0,1$,
then $\langle u_{0},P_{0},v_{0},u_{3},P_{3},v_{3},u_{2},P_{2},v_{2},u_{1},P_{1},v_{1},u_{0}\rangle$ is a fault-free cycle of length $2^{2n}-2f_{v}$ in $BH_{n}$.

Subcase 1.3.1.2. $k=2$.

Let $(u_{0},v_{0})\in C_{0}$ and assume $u_{0}$ is white.
Choose an edge $(u_{2},v_{2})\in E(C_{2})$, assume $u_{2}$ is white.
Let $v_{1}\in V(BH_{n-1}^{1})$ (resp. $u_{1}\in V(BH_{n-1}^{1})$, $v_{3}\in V(BH_{n-1}^{3})$ and $u_{3}\in V(BH_{n-1}^{3})$) be an extra neighbor of $u_{0}$ (resp. $v_{2}$, $u_{2}$ and $v_{0}$).
Since $f_{t}=0$ for $t\in\{1,3\}$, by Lemma~\ref{H-L}, there exists a fault-free path $P_{t}$ of length
$2^{2(n-1)}-1$ between $u_{t}$ and $v_{t}$ in $BH_{n-1}^{t}$. Let $P_{i}=C_{i}-(u_{i},v_{i})$ for $i=0,2$,
then $\langle u_{0},P_{0},v_{0},u_{3},P_{3},v_{3},u_{2},P_{2},v_{2},u_{1},P_{1},v_{1},u_{0}\rangle$ is a desired fault-free cycle of length $2^{2n}-2f_{v}$ in $BH_{n}$.

Subcase 1.3.2. $f_{v}^{0}=n-1$. (The discussion is similar for $f_{v}^{k}=n-1$.)

It implies $f_{v}^{k}=0$ and $f_{e}^{k}=n-1$.
We regard one faulty vertex in $BH_{n-1}^{0}$ as fault-free, $f_{v}'=f_{v}^{0}-1=n-2$.
By the inductive hypothesis, there exists a cycle $C_{0}$ in $BH_{n-1}^{0}$ of length $2^{2(n-1)}-2f_{v}'$, where $f_{v}'=f_{v}^{0}-1=f_{v}-1$. Note that $C_{0}$ contains at most
one faulty vertex. Without loss of generality, assume that $C_{0}$ contains one faulty vertex $u_{0}$ and $u_{0}$ is black. Let $a_{0}$ and $b_{0}$ be two vertices adjacent to $u_{0}$ in $C_{0}$. Let $c_{0}\in C_{0}$ be a neighbor of $b_{0}$.
The desired cycle of $BH_{n}$ can be constructed by the similar discussions as Subcase 1.2.2, the details are omitted.

Case 2. There exists one $i\in \{0,1,2,3\}$ such that $f_{i}=2n-3$.

Without loss of generality, let $i=0$, i.e.,
$f_{0}=2n-3$. It implies that $f_{t}\leq f-f_{0}\leq 1$ for all $t\in \{1,2,3\}$. We consider the following two subcases.

Subcase 2.1.  $f_{v}^{0}\leq n-2$.

Recall that $f_{e}^{0}=f_{0}-f_{v}^{0}\geq 2n-3-(n-2)=n-1\geq 2$ for $n\geq 3$. Choose any faulty edge, say $(u_{0},v_{0})\in E(BH_{n-1}^{0})$,
regard $(u_{0},v_{0})$ as a fault-free edge temporarily, then $f_{0}-1\leq 2n-4$.
By the inductive hypothesis, there exists a cycle $C_{0}$ of length $2^{2(n-1)}-2f_{v}^{0}$ in $BH_{n-1}^{0}$ such that $C_{0}$ contains at most one faulty edge $(u_{0},v_{0})$.
Suppose $(u_{0},v_{0})\in E(C_{0})$ (otherwise choosing an edge in $E(C_{0})$ replace for $(u_{0},v_{0})$). Since there exists at most one faulty element outside $BH_{n-1}^{0}$,
by Lemma~\ref{8-cycle}(3),
$(u_{0},v_{0})$ is contained in a fault-free cycle $C$, say $C=\langle v_{0},u_{0},v_{1},u_{1},v_{2},u_{2},v_{3},u_{3},v_{0}\rangle$, of length $8$
in $BH_{n}$ such that $E(C)\cap E(BH_{n-1}^{i})=\{(u_{i},v_{i})\}$ for each $i\in \{0,1,2,3\}$. Since $f_{t}\leq 1$ for all $t\in \{1,2,3\}$, note that $1\leq n-1$ and $1\leq 2n-3$ for $n\geq 3$,
by Lemma~\ref{vertex-n-1} or Lemma~\ref{edge-2n-3}, there exists a fault-free cycle $C_{t}$ of length $2^{2(n-1)}-2f_{v}^{t}$ passing through $(u_{t},v_{t})$.
Let $P_{i}=C_{i}-(u_{i},v_{i})$ for $i\in \{0,1,2,3\}$,
then $\langle u_{0},P_{0},v_{0},u_{3},P_{3},v_{3},u_{2},P_{2},v_{2},u_{1},P_{1},v_{1},u_{0}\rangle$ is a desired fault-free cycle of length $2^{2n}-2f_{v}$ in $BH_{n}$.

Subcase 2.2.  $f_{v}^{0}=n-1$.

It implies $f_{v}^{j}=0$ and $f_{e}^{j}\leq 1$ for each $j\in \{1,2,3\}$. We regard one faulty vertex in $BH_{n-1}^{0}$ as fault-free temporarily. For the situation where $f_{0}-1\leq 2n-4$
and $f_{v}'=f_{v}^{0}-1=f_{v}-1=n-2$, by the inductive hypothesis, there exists a cycle $C_{0}$ in $BH_{n-1}^{0}$ of length $2^{2(n-1)}-2f_{v}'$,
where $f_{v}'=f_{v}^{0}-1$. Note that $C_{0}$ contains at most
one faulty vertex.

If $C_{0}$ does not contain any faulty vertex, choose an edge, say $(u_{0},v_{0})\in E(C_{0})$
(assume $u_{0}$ is white, $v_{0}$ is black).
By Lemma~\ref{8-cycle}(3) and $f-f_{0}\leq 1$,
$(u_{0},v_{0})$ is contained in a fault-free cycle $C$, say $C=\langle v_{0},u_{0},v_{1},u_{1},v_{2},u_{2},v_{3},\\u_{3},v_{0}\rangle$, of length $8$
in $BH_{n}$ such that $E(C)\cap E(BH_{n-1}^{i})=\{(u_{i},v_{i})\}$ for each $i\in \{0,1,2,3\}$.
By using the method similar to that of Subcase 1.2.2, we can find a desired cycle
of length $2^{2n}-2f_{v}$ in $BH_{n}$. The details are omitted.

Now we assume that $C_{0}$ contains one faulty vertex $u_{0}$ and $u_{0}$ is black (if $u_{0}$ is white, the discussions are similar). Let $a_{0}$
and $b_{0}$ be two vertices adjacent to $u_{0}$ in $C_{0}$. Let $c_{0}\in C_{0}$ be a neighbor of $b_{0}$.
Recall that each vertex has two extra neighbors and $f_{c}\leq 1$, we can choose the extra neighbors $c_{1}\in V(BH_{n-1}^{1})$ and $a_{3}\in V(BH_{n-1}^{3})$ of $a_{0}$ and $c_{0}$, respectively,
such that both $(a_{0},c_{1})$ and $(c_{0},a_{3})$ are fault-free.
The desired cycle of $BH_{n}$ can be constructed by the similar discussions as Subcase 1.2.2, the details are omitted.

Case 3. There exists one $i\in \{0,1,2,3\}$ such that $f_{i}=2n-2$.

Without loss of generality, let $i=0$, i.e. $f_{0}=f=2n-2$. It implies that $f_{j}=0$ for all $j\in \{1,2,3\}$ and $f_{c}=0$.

Since $f_{v}\geq 1$, we can choose a faulty vertex, say $w$ (assume $w$ is black, the discussion is similar for $w$ being white). Note that $f_{v}\leq n-1$, $f_{e}=f-f_{v}\geq 2n-2-(n-1)=n-1\geq 1$ for $n\geq 3$,
let $(u_{0},v_{0})$ (assume $v_{0}$ is black) be a faulty edge. Image $w$ and $(u_{0},v_{0})$ as fault-free temporarily.
Since $f_{0}-2\leq 2n-4$
and $f_{v}'=f_{v}^{0}-1\leq n-1-1=n-2$, by the inductive hypothesis, there exists a cycle $C_{0}$ in $BH_{n-1}^{0}$ of length $2^{2(n-1)}-2f_{v}'$,
where $f_{v}'=f_{v}^{0}-1=f_{v}-1$. Note that $C_{0}$ contains at most
one faulty vertex $w$ and one faulty edge $(u_{0},v_{0})$. We consider the following four subcases.

Subcase 3.1. $w\notin C_{0}$ and $(u_{0},v_{0})\notin E(C_{0})$.

Choose any edge, say $(a_{0},b_{0})\in E(C_{0})$ (assume $a_{0}$ is white, $b_{0}$ is black).
By Lemma~\ref{8-cycle}(1) and $f_{c}=0$,
$(a_{0},b_{0})$ is contained in a fault-free cycle $C$, say $C=\langle b_{0},a_{0},b_{1},a_{1},b_{2},a_{2},b_{3},a_{3},a_{0}\rangle$, of length $8$
in $BH_{n}$ such that $E(C)\cap E(BH_{n-1}^{i})=\{(a_{i},b_{i})\}$ for each $i\in \{0,1,2,3\}$.
Since for any $i\in \{1,2\}$, $f_{i}=0$, by Lemma~\ref{H-L},
there exists a fault-free path $P_{i}$ of length $2^{2(n-1)}-1$ between $a_{i}$ and $b_{i}$ in $BH_{n-1}^{i}$ for $i=1,2$.
Since $f_{3}=0$, by Lemma~\ref{bipan}, there exists a fault-free path $P_{3}$ of length $2^{2(n-1)}-3$ between $a_{3}$ and $b_{3}$ in $BH_{n-1}^{i}$.
Let $P_{0}=C_{0}-(a_{0},b_{0})$.
A desired fault-free cycle of $2^{2n}-2f_{v}$ in $BH_{n}$ can be constructed
as $\langle a_{0},b_{1},P_{1},a_{1},b_{2},P_{2},a_{2},b_{3},P_{3},a_{3},b_{0},P_{0},a_{0}\rangle$.

Subcase 3.2. $w\notin C_{0}$ and $(u_{0},v_{0})\in E(C_{0})$.

By Lemma~\ref{8-cycle}(1) and $f_{c}=0$,
$(u_{0},v_{0})$ is contained in a cycle $C$, say $C=\langle v_{0},u_{0},v_{1},u_{1},v_{2},\\u_{2},v_{3},u_{3},v_{0}\rangle$, of length $8$
in $BH_{n}$ such that $E(C)\cap E(BH_{n-1}^{i})=\{(u_{i},v_{i})\}$ and $(u_{i},v_{i+1})$ is fault-free for each $i\in \{0,1,2,3\}$.
A desired cycle can be constructed by the similar discussion as Subcase 3.1, the details are omitted.

\begin{figure}

\begin{minipage}[t]{0.5\linewidth}
\centering
\includegraphics[height=5cm,width=5cm]{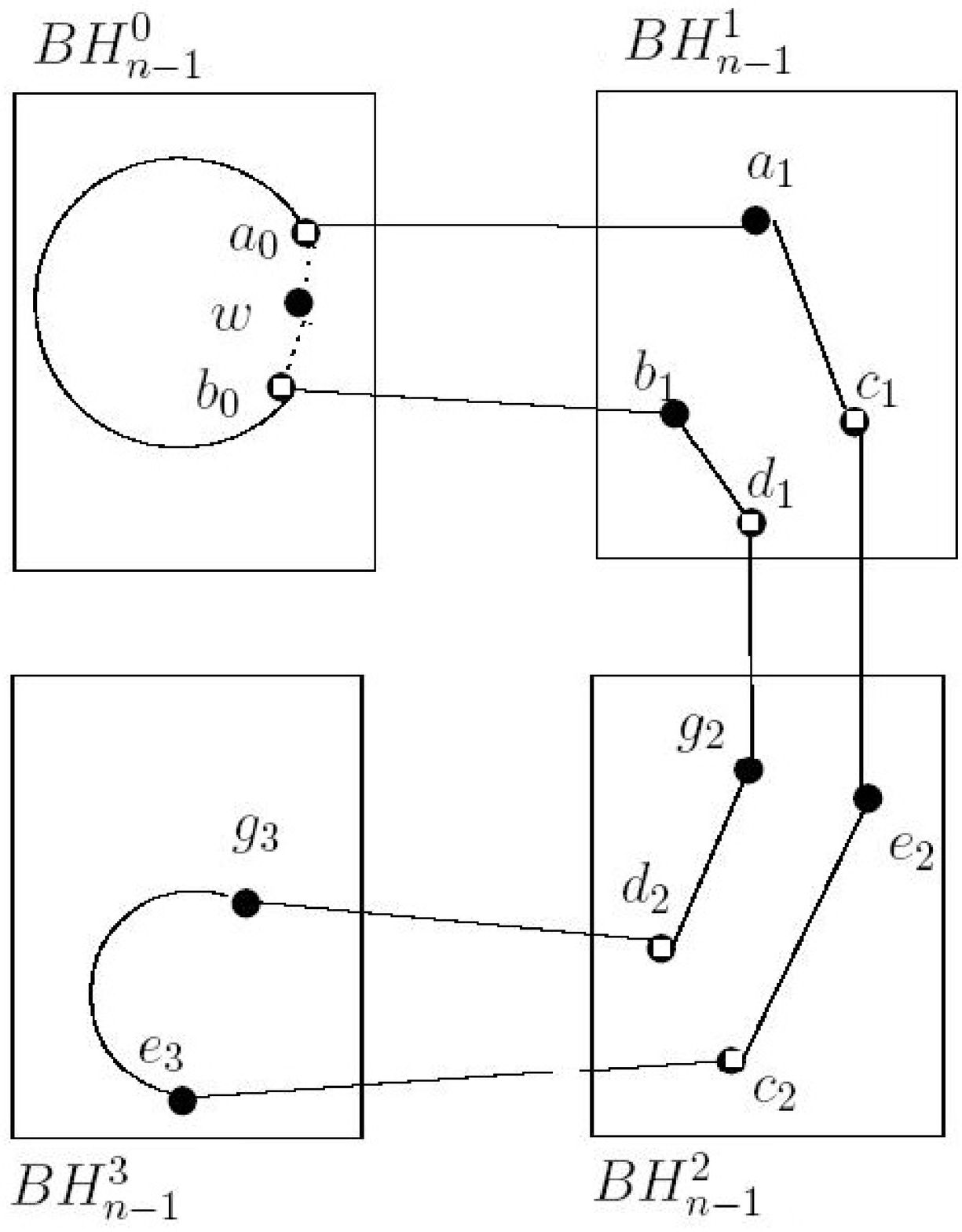}
\vskip-0.3cm
\caption{Illustration of 3.3;}
\label{3.3}
\end{minipage}
\begin{minipage}[t]{0.5\linewidth}
\centering
\includegraphics[height=5cm,width=5cm]{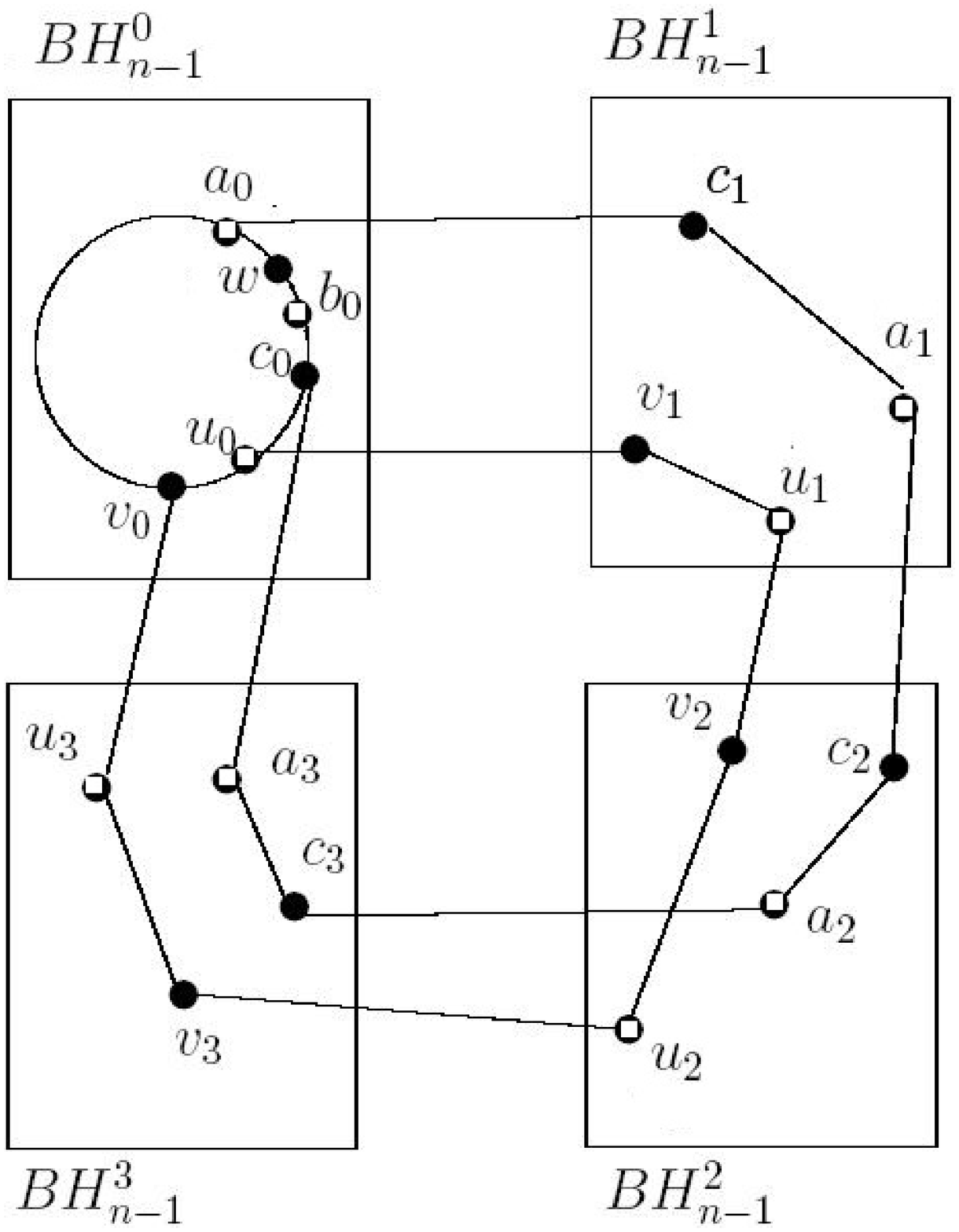}
\vskip-0.3cm
\caption{Illustration of 3.4.2}
\label{3.4.2}
\end{minipage}

\end{figure}

Subcase 3.3. $w\in C_{0}$ and $(u_{0},v_{0})\notin E(C_{0})$.

Let $a_{0}$ and $b_{0}$ be two vertices adjacent to $w$ in $C_{0}$. Since each vertex has two extra neighbors,
let $a_{1}\in V(BH_{n-1}^{1})$, $b_{1}\in V(BH_{n-1}^{1})$ be an extra neighbor of $a_{0}$ and $b_{0}$, respectively, such that $a_{1}\neq b_{1}$.
Choose two distinct white vertices, say $c_{i}$, $d_{i}$, respectively, in $BH_{n-1}^{i}$ for $i=1,2$.
Let one of the extra neighbors of $c_{i}$, $d_{i}$ be $e_{i+1}$, $g_{i+1}$ in $BH_{n-1}^{i+1}$, respectively, for $i=1,2$.
See Figure~\ref{3.3}. By Lemma~\ref{two-node}, there exist two vertex-disjoint paths $P_{1}[a_{1},c_{1}]$ and $R_{1}[b_{1},d_{1}]$ (resp. $P_{2}[c_{2},e_{2}]$ and $R_{2}[d_{2},g_{2}]$ )
in $BH_{n-1}^{1}$ (resp. $BH_{n-1}^{2}$) such that $V(P_{1}[a_{1},c_{1}])\cup V(R_{1}[b_{1},d_{1}])=V(BH_{n-1}^{1})$ (resp. $V(P_{2}[c_{2},e_{2}])\cup V(R_{2}[d_{2},g_{2}])=V(BH_{n-1}^{2})$).
By Lemma~\ref{laceable}, there exists a path $P_{3}$ in $V(BH_{n-1}^{3})$ of length $2^{2(n-1)}-2$ between $e_{3}$ and $g_{3}$.
Since $f_{c}=0$, the crossing edges $(c_{1},e_{2})$, $(d_{1},g_{2})$, $(c_{2},e_{3})$, $(d_{2},g_{3})$ are fault-free. Let $P_{0}=C_{0}-\{w\}$.
Thus, a fault-free cycle of length $2^{2n}-2f_{v}^{0}=2^{2n}-2f_{v}$ in $BH_{n}$ is $\langle a_{0},a_{1},P_{1}[a_{1},c_{1}],c_{1},e_{2},P_{2}[e_{2},c_{2}],c_{2},e_{3},P_{3},g_{3},d_{2},
R_{2}[d_{2},g_{2}],g_{2},d_{1},R_{1}[d_{1},b_{1}],\\b_{1},b_{0},P_{0},a_{0}\rangle$.

Subcase 3.4. $w\in C_{0}$ and $(u_{0},v_{0})\in E(C_{0})$.

Let $a_{0}$ and $b_{0}$ be two vertices adjacent to $w$ in $C_{0}$. If $v_{0}=w$, then $u_{0}=a_{0}$ or $u_{0}=b_{0}$. Without loss of generality, assume that $u_{0}=b_{0}$. The discussion is the same as that in Subcase 3.3.

Now we assume $v_{0}\neq w$ in the following. Notice that $u_{0}$ may equal to $a_{0}$ or $b_{0}$.

Subcase 3.4.1. $u_{0}=b_{0}$ (If $u_{0}=a_{0}$, the discussion is similar).

Assume that $C_{0}=\langle a_{0},w,u_{0},v_{0},P_{0}[v_{0},a_{0}],a_{0}\rangle$.
For $i\in \{1,2,3\}$, let $u_{i}\in V(BH_{n-1}^{i})$ be a white vertex such that $(u_{3},v_{0})\in E(BH_{n})$.
For $j\in \{1,2,3\}$, let $v_{j}\in V(BH_{n-1}^{j})$ be a black vertex such that $(a_{0},v_{1}),(u_{1},v_{2}),(u_{2},\\v_{3})\in E(BH_{n})$.
Since $f_{i}=0$ for $i\in\{1,2,3\}$, by Lemma~\ref{H-L}, there exists a $u_{i}-v_{i}$ path $P_{i}$ of length $2^{2(n-1)}-1$ in $BH_{n-1}^{i}$.
Thus, $\langle a_{0},v_{1},P_{1},u_{1},v_{2},P_{2},u_{2},\\v_{3},P_{3},u_{3},v_{0},P_{0}[v_{0},a_{0}],a_{0}\rangle$ is a desired fault-free cycle of length $2^{2n}-2f_{v}$ in $BH_{n}$.

Subcase 3.4.2. $u_{0}\neq \{a_{0},b_{0}\}$.

Let $c_{1}\in V(BH_{n-1}^{1})$ (resp. $u_{3}\in V(BH_{n-1}^{3})$) be an extra neighbor of $a_{0}$ (resp. $v_{0}$).
Let $c_{0}\in C_{0}$ be a neighbor of $b_{0}$ such that $c_{0}\neq v_{0}$ (otherwise let $c_{0}\in C_{0}$ be a neighbor of $a_{0}$, the discussion is similar). Thus $C_{0}$ is the type as $C_{0}=\langle a_{0},w,b_{0},c_{0},P_{0}[c_{0},u_{0}],u_{0},v_{0},P_{0}[v_{0},a_{0}],a_{0}\rangle$ (to see Figure~\ref{3.4.2}).
Since each vertex has two extra neighbors, let $v_{1}\in V(BH_{n-1}^{1})$ (resp. $a_{3}\in V(BH_{n-1}^{3})$) be an extra neighbor of $u_{0}$ (resp. $c_{0}$) such that $v_{1}\neq c_{1}$ (resp. $a_{3}\neq u_{3}$).
For $i\in \{1,2\}$, let $u_{i},a_{i}\in V(BH_{n-1}^{i})$ be two distinct white vertices. For $j\in \{2,3\}$, let $v_{j},c_{j}\in V(BH_{n-1}^{j})$ be two distinct black vertices such that
$(u_{j-1},v_{j})\in E(BH_{n})$ and $(a_{j-1},c_{j})\in E(BH_{n})$.
Since $f_{i}=0$ for $i\in \{1,2,3\}$, by Lemma~\ref{two-node},
there exist two vertex-disjoint paths $P_{i}[v_{i},u_{i}]$ and $R_{i}[c_{i},a_{i}]$ in $BH_{n-1}^{i}$ such that $V(P_{i}[v_{i},u_{i}])\cup V(R_{i}[c_{i},a_{i}])=V(BH_{n-1}^{i})$.
Thus, a fault-free cycle of length $2^{2n}-2f_{v}$ in $BH_{n}$ is
$\langle a_{0},c_{1},R_{1}[c_{1},a_{1}],a_{1},c_{2},R_{2}[c_{2},a_{2}],a_{2},c_{3},R_{3}[c_{3},a_{3}],a_{3},c_{0},P_{0}[c_{0},u_{0}],\\
u_{0},v_{1},P_{1}[v_{1},u_{1}],u_{1},v_{2},P_{2}[v_{2},u_{2}],u_{2},v_{3},P_{3}[v_{3},u_{3}],u_{3},v_{0},P_{0}[v_{0},a_{0}],a_{0}\rangle$.

\medskip

By the above cases, the proof is complete. \hfill\qed

\section{Conclusion}

In this paper, we obtained that $BH_{n}$ with $|F_{v}|+|F_{e}|\leq 2n-2$ and $|F_{v}|\leq n-1$ has a fault-free cycle of length $2^{2n}-2|F_{v}|$.
If $|F_{v}|=0$, $|F_{e}|\leq 2n-2$, a fault-free cycle of length $2^{2n}-2|F_{v}|=2^{2n}$ is a Hamiltonian cycle.
Since $BH_{n}$ is a bipartite graph with two
partite sets of equal size, the cycle is the longest in the worst-case. Furthermore, since the edge
connectivity of $BH_{n}$ is $2n$, the $BH_{n}$ cannot tolerate $2n-1$ faulty edges. Hence,
our result is optimal in terms of faulty edges.

\medskip

\section*{Acknowledgments}

This work was supported by the National Natural Science Foundation of China (No.11371052, No.11571035, No.11271012 and No. 61572005).


\end{spacing}

\begin{thebibliography}{99}



\bibitem{A}
B. Alspach, D. Hare, Edge-pancyclic block-intersection graphs, Discrete Math. 97 (1991) 17--24.

\bibitem{B}
J.A. Bondy, Pancyclic graphs, J. Combin. Theory 11 (1971) 80--84.

\bibitem{C}
D. Cheng, R.-X. Hao, Various cycles embedding in faulty balanced hypercubes, Inform. Sci. 297 (2015) 140--153.

\bibitem{C1}
D. Cheng, R.-X. Hao, Y.-Q. Feng, Vertex-fault-tolerant cycles embedding in balanced hypercubes, Inform. Sci. 288 (2014) 449--461.

\bibitem{C2}
D. Cheng, R.-X. Hao, Y.-Q. Feng, Two node-disjoint paths in balanced hypercubes, Appl. Math. Comput. 242 (2014) 127--142.

\bibitem{D1}
Z.-Z. Du, J.-M. Xu, A note on cycle embedding in hypercubes with faulty vertices, Inform. Proces. Lett. 111 (2011) 557--560.


\bibitem{F}
J.-S. Fu, Longest fault-free paths in hypercubes
with vertex faults, Inform. Sci. 176 (2006) 759--771.


\bibitem{H}
R.-X. Hao, R. Zhang, Y.-Q. Feng, J.-X. Zhou, Hamiltonian cycle embedding for fault tolerance in balanced
hypercubes, Appl. Math. Comput. 244 (2014) 447--456.


\bibitem{H2}
S.-Y. Hsieh, Fault-tolerant cycle embedding in the hypercube
with more both faulty vertices and faulty edges, Parall. Comput. 32 (2006) 84--91.



\bibitem{H3}
K. Huang, J. Wu, Balanced hypercubes, in: Proc. of the 1992, Int. Conf. Parallel Process. 3 (1992) 153--159.

\bibitem{H4}
K. Huang, J. Wu, Fault-tolerant resource placement in balanced hypercubes, Inform. Sci. 99 (1997) 159--172.

\bibitem{H5}
K. Huang, J. Wu, Area efficient layout of balanced hypercubes, Int. J. High Speed Electron. Syst. 6 (1995) 631--646.



\bibitem{L1}
H. L\"{u}, X. Li, H. Zhang, Matching preclusion for balanced hypercubes,
Theoret. Comput. Sci. 465 (2012) 10--20.

\bibitem{L2}
H. L\"{u},  H. Zhang, Hyper-Hamiltonian laceability of balanced hypercubes, J. Supercomput. 68 (1) (2014) 302--314.



\bibitem{T}
C.-H. Tsai, Fault-tolerant cycles embedded in hypercubes with mixed link and node failures, Appl. Math. Lett. 21 (2008) 855--860.


\bibitem{W}
J. Wu, K. Huang, The balanced hypercubes: a cube-based system for fault-tolerant applications, IEEE Trans. Comput. 46 (4) (1997) 484--490.

\bibitem{X0}
J.-M. Xu, M. Ma, Survey on path and cycle embedding in some networks, Front. Math. China 4 (2009) 217--252.

\bibitem{X}
M. Xu, X.-D. Hu, J.-M. Xu, Edge-pancyclicity and Hamiltonian laceability
of the balanced hypercubes, Appl. Math. Comput. 189 (2007)
1393--1401.

\bibitem{Y1}
M.-C. Yang, Bipanconnectivity of balanced hypercubes, Comput. Math. Appl. 60 (2010) 1859--1867.

\bibitem{Y2}
M.-C. Yang, Super connectivity of balanced hypercubes, Appl. Math. Comput. 219 (2012) 970--975.

\bibitem{Y3}
M.-C. Yang, Conditional diagnosability of balanced hypercubes under the PMC model, Inform. Sci. 222 (2013) 754--760.


\bibitem{Ym}
M.-C. Yang, J.J.M. Tan, L.H. Hsu, Hamiltonian circuit and linear array embeddings in faulty $k$-ary $n$-cubes, J. Parall. Distrib. Comput. 67 (4) (2007) 362--368.

\bibitem{Y4}
M.-C. Yang, M.-H. Yang, Reliability analysis of balanced hypercubes,
IEEE COMCOMAP, (2012) 376--379.

\bibitem{Y5}
W.-H. Yang, H.Z. Li, W.H. He, Edge-fault-tolerant bipancyclicity of Cayley graphs generated
by transposition-generating trees, Inter. J. Comput. Math. 92(2015) 1345--1352.

\bibitem{Z}
J.-X. Zhou, Z.-L. Wu, S.-C. Yang, K.-W. Yuan, Symmetric property
and reliability of balanced hypercube, IEEE Trans. Comput. 64 (2015) 876--881.

\bibitem{Z1}
Q. Zhou, D. Chen, H. L\"{u}, Fault-tolerant Hamiltonian laceability of balanced hypercubes, Inform. Sci. 300 (2015) 20--27.



\end{thebibliography}
\end{document}